\documentclass[final,12pt,3p]{elsarticle}

\usepackage{amssymb}
\usepackage{amsthm}
\usepackage{caption}
\usepackage{subcaption}
\usepackage{mathtools}
\usepackage{hyperref}
\usepackage[utf8]{inputenc}
\usepackage{graphicx}
\usepackage{lineno}
\usepackage[procnames]{listings}
\usepackage[bottom,hang]{footmisc}
\usepackage{xcolor}
\usepackage{algorithm,algpseudocode}
\algnewcommand{\Inputs}[1]{%
  \State \textbf{Inputs:}
  \Statex \hspace*{\algorithmicindent}\parbox[t]{.8\linewidth}{\raggedright #1}}
\algnewcommand{\Initialize}[1]{%
  \State \textbf{Initialize:}
  \Statex \hspace*{\algorithmicindent}\parbox[t]{.8\linewidth}{\raggedright #1}}
\algnewcommand{\Output}[1]{%
  \State \textbf{Output:}
  \Statex \hspace*{\algorithmicindent}\parbox[t]{.8\linewidth}{\raggedright #1}}
\definecolor{keywords}{RGB}{0,0,255}
\definecolor{comments}{RGB}{80,80,150}
\definecolor{blue}{RGB}{0,0,255}
\definecolor{black}{RGB}{0,0,0}

\journal{~}

\begin{document}

\begin{frontmatter}

\title{Cost-Optimal Operation of Energy Storage Units: \\ Benefits of a Problem-Specific Approach}

\author[label1]{Lars Siemer\corref{cor1}}
\ead{lars.siemer@next-energy.de}
\address[label1]{NEXT ENERGY -- EWE Research Centre for Energy Technology at the University of Oldenburg, Carl-von-Ossietzky-Str.~15, 26129 Oldenburg, Germany}

\address[label2]{Department of Mathematics, University of Oldenburg, 26129 Oldenburg, Germany}

\cortext[cor1]{Corresponding author}

\author[label2]{Frank Sch{\"o}pfer}

\author[label1]{David Kleinhans}

\begin{abstract}
The integration of large shares of electricity produced by non-dispatchable Renewable Energy Sources (RES) leads to an increasingly volatile energy generation side, with temporary local overproduction. The application of energy storage units has the potential to use this excess electricity from RES efficiently and to prevent curtailment. The objective of this work is to calculate cost-optimal charging strategies for energy storage units used as buffers. For this purpose, a new mathematical optimization method is presented that is applicable to general storage-related problems. Due to a tremendous gain in efficiency of this method compared with standard solvers and proven optimality, calculations of complex problems as well as a high-resolution sensitivity analysis of multiple system combinations are feasible within a very short time. As an example technology, Power-to-Heat converters used in combination with thermal storage units are investigated in detail and optimal system configurations, including storage units with and without energy losses, are calculated and evaluated. The benefits of a problem-specific approach are demonstrated by the mathematical simplicity of our approach as well as the general applicability of the proposed method.
\end{abstract}
\end{frontmatter}

\section{Introduction}
\label{sec:sec1}
Electricity generation from solar and wind resources is volatile and has limited adaptability to changes in electricity demand \cite{georgilakis2008technical}. For a stable operation of the power system there must be a consistent balance between supply and demand. One option for achieving this balance is the use of energy storage units which can contribute substantially to the expansion of Renewable Energy Sources (RES) and their integration into existing grids \cite{heide2010seasonal, rasmussen2012storage, weitemeyer2015integration}. 

Another possible way in which large shares of renewable energy can be integrated and electricity from RES can be used efficiently, especially in times of overproduction, is to replace fossil fuels as the main energy source for heat by converting electricity into heat (Power-to-Heat, PtH) -- resulting in a lower overall primary energy consumption and lower CO\(_2\) emissions \cite{potenzial2014, sterner2014energiespeicher}. If combined with thermal storage units, PtH storage systems are a highly flexible option for uncoupling conversion and utilization \cite{potenzial2014, denholm2010role}. These heat storage units can be charged at considerably low losses during periods with high shares of excess electricity from RES and provide heat in times of low feed-in \cite{pedersen2011}. A method for storing heat energy used in Germany was implemented decades ago by night storage heating systems in private households. They operated with low overall efficiencies and subsidized electricity prices at fixed schedules \cite{eiselt2012dezentrale}. Today, private households are still suitable for the deployment of PtH storage systems due to the high share of primary energy demand utilized for heating and hot water supply, the overall changes in the energy system and the characteristics of today's heat storage units \cite{stat2014energie}.

If we assume that the electricity price represents the fluctuating availability of energy, these heat storage systems need to be operated in a cost-optimal manner. The optimality of a charging strategy for a PtH storage unit is determined by the overall electricity acquisition costs. The task of identifying optimal strategies is closely related to the field of mathematical optimization and can be described as a minimization problem. Standard solvers can be applied to calculate solutions, albeit with tremendously high overheads in computational time. Hence, the investigation of complex problems may not be possible within a reasonable time frame.

The specific structure of the minimization problem, based on a mathematical model of the storage units, allows the development of a new optimization method which to our knowledge has not yet been implemented. This work introduces an innovative optimization algorithm and presents proof of optimality. The solution of storage-related optimization problems can now be calculated in a fraction of the time required by standard algorithms. As a particular example, we focus on PtH storage units installed in private households. For these systems, cost-optimal operating strategies including storage units with thermal energy losses are described and used for an iterative determination of optimal system designs.

Beyond this particular field of application, we illustrate the potential benefits of a problem-specific optimization approach in terms of mathematical simplicity and computational gain as compared with standard solvers.


\section{Formulation of the optimization problem}
\label{sec:sec2}

The construction of mathematical energy storage models and the formulation of the corresponding optimization problems, with and without constant as well as non-constant energy losses, are described without units in the following paragraphs. We investigate an electrical charging system with an attached storage unit as a buffer and consider the parameters of electric power consumption and storage capacity. The system needs to be connected to a virtual electricity grid which provides an altering price signal reflecting the availability of electricity. A further precondition is that the energy demand is covered exclusively by the system. 

\subsection{Formulation without energy losses}
\label{subsec:no losses}

Now, we formulate the mathematical model, without energy losses of the storage unit over time or losses during the conversion process. We discretize the considered period of time into $n$ intervals of equal size. The electricity prices are denoted by \(c_1,c_2, \dots , c_n\) and the energy demands are represented by \(d_1, d_2, \dots , d_n \). Furthermore, we define the amount of energy used to charge the storage unit for each time-interval by \(x_1,x_2, \dots , x_n\) and assume that the prices and the demands are known. For technical reasons, the values of $x$ are constrained. Each charge value cannot be negative (no discharge to the electricity grid or re-electrification) and has an upper bound \(C>0\), which implies:
\[0 \le x_i \le C\,, \ \forall \, i=1, \dots ,n\,.\]
To cover at least the energy demand for each time-interval, additional constraints are:
\[\sum_{j=1}^{i}{d_j} \le \sum_{j=1}^{i}{x_j}\,, \ \forall \, i=1, \dots ,n\,.\]
The storage level for each time-interval is defined as the difference between the quantity of charges and the quantity of demands up to this time-interval. Due to design-related restrictions, the storage level is bounded above by a maximum storage capacity value \(S>0\). Hence, we have:
\[\sum_{j=1}^{i}{\left(x_j-d_j\right)} \le S\,, \ \forall \, i=1, \dots ,n\,.\]
In order to minimize the total costs to cover at least the energy demand over a period of time of size $n$, we have to solve the following optimization problem:
\begin{align}
\min \quad & \sum_{i=1}^{n}{c_i \cdot x_i}   \nonumber \\
\mbox{subject to} \quad & 0 \le x_i \le C\,, \ \mbox{$\forall \, i=1, \dots ,n$} \label{eq:LP} \tag{LP} \\
\mbox{and}\quad  & \sum_{j=1}^{i}{d_j} \le \sum_{j=1}^{i}{x_j} \le S+\sum_{j=1}^{i}{d_j}\,, \ \mbox{$\forall \, i=1, \dots ,n$} \,.  \nonumber
\end{align}
This problem is a special instance of a so-called \emph{linear program} and can be solved using standard algorithms for general linear programs, such as the simplex method or interior point methods \cite{nocedal2006numerical}.

The problem~\eqref{eq:LP} is also a special case of the problem~\eqref{eq:P} which is described in~\ref{app:algorithm}. Contrary to the methods mentioned above, we utilize the special structure of~\eqref{eq:P} and hence also of the related problem~\eqref{eq:LP} to develop a new algorithm. The basic idea of the new algorithm is to charge the storage unit during periods when the acquisition prices are low in order to avoid further purchases at times when the prices are higher. In addition to the price levels, the algorithm also takes into account the demand, the storage level and the maximum charge power for each time-interval. Therefore, the storage units are charged as much as possible at times of negative acquisition costs and as much as required, if the price is non-negative.

The new algorithm to solve the problem~\eqref{eq:P} is discussed in detail in~\ref{app:algorithm} and also the pseudo-code is presented (see p.\ \pageref{alg}). Below, we present the pseudo-code of the new algorithm exemplary fitted to the problem \eqref{eq:LP}, if we set \(a_i:=d_1+ \ldots + d_i\) and \(b_i:=S+\left(d_1+ \ldots + d_i \right)\) for all \(i=1,\dots ,n\). If we further define the permutation \(\sigma\) as described in~\ref{app:algorithm}, corresponding to the increasing prices by \(c_{\sigma(1)} \le \ldots \le c_{\sigma(n)}\), the pseudo-code of the new algorithm fitted to the problem~\eqref{eq:LP} is given by:

\begin{algorithm}[H]
    \begin{algorithmic}
        \vspace{12pt}
        \Inputs{$a_i,b_i,c_i$ and the permutation $\sigma$}
        \Output{A solution \(x\) of the problem~\eqref{eq:LP}}
        \For{\(k=1\) to \(n\)}
	        \State \(M_1 \gets \max\limits_{i<\sigma(k)}\{0,a_i\}\)
            \State \(M_2 \gets \max\limits_{i \ge \sigma(k)}\{0,a_i\}\)
            \State \(m \gets \min\limits_{i \ge \sigma(k)} \{b_i\}\)
	        \If{\(c_{\sigma(k)} \ge 0\)}
				\State $x_{\sigma(k)} \gets \min\big\{\max\{0,M_2-M_1\},\min\{C, m-M_1\}\big\}$
			\Else
				\State $x_{\sigma(k)} \gets \min\{C, m-M_1\}$
			\EndIf
			\For{$i=\sigma(k),\ldots,n$}
		    \State $a_i \gets a_i-x_{\sigma(k)}$
		    \State $b_i \gets b_i-x_{\sigma(k)}$
	    \EndFor
    \EndFor
    \vspace{12pt}
  \end{algorithmic}
\end{algorithm}
\vspace{12pt}

The optimality of the new algorithm on page~\pageref{app:algorithm} is proven (in~\ref{app:algorithm}) and is one key element of this work. Furthermore, the source code of an implementation in Python is included in~\ref{app:python}. This method solves problems even for large $n$ in a fraction of time required by standard solvers, because at most \(n^2 +3n\) floating point operations and \(3/2n^2+5/2n\) comparisons to compute a solution (see:~\ref{app:algorithm}) are required.

To demonstrate the efficiency of the new algorithm, we compare the runtimes between the algorithm and the common solver \emph{linprog} in the following. For this comparison, a straightforward implementation of the new algorithm in Python (cf.\ \ref{app:python}) and the \emph{linprog} implementation as available in MATLAB 2015b were used. The calculations were performed on a desktop computer\footnote{Details:
Intel$^{\textsuperscript{\textregistered}}$ Core$^{\text{TM}}$ i7-930 Processor (2.80 GHz), 12GB RAM, MATLAB 2015b, Python 2.7.}, based on input data for the problem \eqref{eq:P} as described in section~\ref{subsec:subsec1}. The results are presented in table~\ref{table:runtime} and emphasize the efficiency of the new algorithm and its beneficial scaling with the size of the problem. This efficiency among others allows to perform detailed sensitivity analyses and solve highly complex optimization problems, where problem~\eqref{eq:P} has to be solved thousands of times as a sub-problem (cf.\ section~\ref{sec:sec3}). For instance, the calculations performed in section~\ref{subsec:subsec2} would have taken more than 60 days by using the solver \emph{linprog} compared with about 25 minutes for the newly proposed algorithm. 

\renewcommand{\arraystretch}{1.3}
\begin{table}
\centering
\begin{tabular}{l|c c c c c}
\hline
n & 100 & 500 & 1000 & 2500 & 5000 \\
\hline
\hline
Runtime algorithm [s] & 0.006 & 0.03 & 0.07 & 0.18 & 0.4 \\
Runtime \emph{linprog} [s] &  1.008  &  2.41 & 18.16  &  232.31 & 1308.9 \\
\hline
\end{tabular}
\caption{Calculation times of the new algorithm compared with \emph{linprog} (MATLAB) of the same problem \eqref{eq:LP}, with respect to the problem size $n$. In addition to the general efficiency of the new algorithm, its beneficial scaling behavior (runtime increases by a factor $65$ if the problem size is increased from $100$ to $5000$, whereas the scaling factor for \emph{linprog} is about 1300) becomes evident.}
\label{table:runtime}
\end{table}

To include energy losses of the storage unit or during the conversion process, which can be constant or not-constant, the constraints of~\eqref{eq:LP} have to be transformed into the shape of the general problem~\eqref{eq:P}, before the new algorithm can be applied. These modifications are presented below for the two different types which can also be combined.

\subsection{Constant energy losses}
\label{subsec:constant lossses}

To include constant energy losses \(l_1, \dots , l_n\) for each time-interval, the constraints in~\eqref{eq:LP} have to be modified. These losses can simply be added to the energy demand for each time-interval, which leads to the following problem, with the same structure as \eqref{eq:P}:
\begin{align}
\min \quad & \sum_{i=1}^{n}{c_i \cdot x_i}   \nonumber \\
\mbox{subject to} \quad & 0 \le x_i \le C\,, \ \mbox{$\forall \, i=1, \dots ,n$} \label{eq:LP constant}\\
\mbox{and}\quad  & \sum_{j=1}^{i}{d_j+l_i} \le \sum_{j=1}^{i}{x_j} \le S+\sum_{j=1}^{i}{d_j+l_i}\,, \ \mbox{$\forall \, i=1, \dots ,n$} \,.  \nonumber
\end{align}
\subsection{Non-constant energy losses}
\label{subsec:variable losses}

A slightly more complicated modification of~\eqref{eq:LP} involves the integration of non-constant energy losses, which is described in~\ref{subsec:appsubsec} in detail. This represents for instance losses during the charging process or energy losses of the storage unit over time. The modified problem is:
\begin{align}
\min \quad & \sum_{i=1}^{n}{c_i \cdot x_i}   \nonumber \\
\mbox{subject to} \quad & 0 \le x_i \le C\,, \ \mbox{$\forall \, i=1, \dots ,n$} \label{eq:LP variable}\\
\mbox{and}\quad  & \sum_{j=1}^{i}{q^{i-j} \cdot d_j} \le \sum_{j=1}^{i}{q^{i-j} \cdot x_j} \le S+\sum_{j=1}^{i}{q^{i-j} \cdot d_j}\,, \ \mbox{$\forall \, i=1, \dots ,n$} \,,  \nonumber
\end{align}
with an energy loss factor \(q>0\). This problem can be transformed into an equivalent problem of the form \eqref{eq:P}, which is described in \ref{subsec:appsubsec} in detail.

Note that the minimization problems \eqref{eq:LP}, \eqref{eq:LP constant} and \eqref{eq:LP variable} are special instances of the problem~\eqref{eq:P} (page \pageref{app:algorithm}). Because the new algorithm (page~\pageref{alg}) solves the general problem~\eqref{eq:P}, it also solves the above mentioned problems. Thus, this algorithm will be used for all discussions and calculations in the following investigations.


\section{Applications to PtH storage units and numerical results}
\label{sec:sec3}

The technical realization of a PtH conversion unit can be achieved for example by an electric heat pump system or an electrical heating element, using the concept of an immersion heater. Common conversion systems offer almost 100~$\%$ efficiency for a simple heating element or even more for a heat pump system, with a performance factor of three to four \cite{von1981heat, perko2011calculation}. Therefore, we do not include energy losses of the PtH system during the conversion process. The thermal storage unit, which is assumed to be attached to the PtH system, can be realized for example by one of these three storage groups: sensitive, latent, or thermochemical \cite{sterner2014energiespeicher, dincer2002thermal, schroder1981latent}. For simplicity, the relevant parameters are the electric power consumption for the converter, the different storage capacities and the storage unit efficiencies.

This section presents a procedure to determine an example of a heat demand profile and calculations of cost-optimal charge curves via the new algorithm for storage units with and without energy losses. Furthermore, we present an approach to estimate cost-optimal combinations of heating systems and connected storage units, with different efficiencies of the storage unit.  

To perform a detailed analysis of the demand for thermal energy in private households, more accurate values than the specific heat demand for one year are required. A common option for calculating these load profiles is to use the VDI Guideline 4655 ``reference load profiles of single and multi-family homes for the use of CHP plants''. The guideline includes the following influence variables to determine such profiles: type of household, number of residents, annual energy demand for space heating and hot water supply as well as the climate zone. 

The majority of households in Germany in 2011 were single and double-family homes, with an average living area of 116.6~m$^2$ \cite{zensus2011b}. The number of households with two residents ranked second with more than 33~$\%$ in the same year and with the highest rate of increase for the future \cite{zensus2011a, entwicklung2010}. In this work, we study a single-family building with two residents. The specific heat demand of these households was on average around 150~kWh/(m$^2 \cdot $a) in 2011, but displayed a decreasing tendency when compare to previous years \cite{bigalke2012dena}. With respect to the EnEV 2009 and the sustained decreasing specific heat demands, a value of 100~kWh/(m$^2 \cdot $a) is considered \cite{bigalke2012dena, EnEV2009}. The guideline recommends an annual heat demand for hot water supply of 500~kWh/(person$\cdot$a), which depends considerably on the number of residents of the observed household and their habits. However, it can be assumed that this value underestimates today's consumption and the adoption of a significantly higher number for real averages approximates this value more closely. In this work, we assume a value of 750~kWh/(person$\cdot$a).

\begin{figure}
\centering
\includegraphics[width=\textwidth]{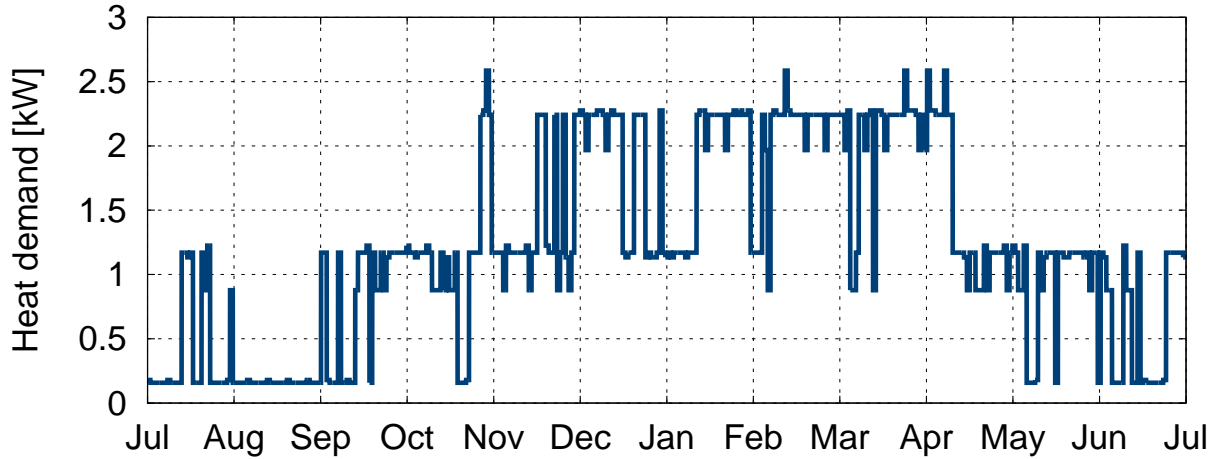}
\caption{Hourly heat demand (space and water heating) profile for a single-family household with two persons, living space of 100~m$^2$, specific heat demand of 100~kWh/(m$^2 \cdot $a) and annual heat demand of 10741.15~kWh/a, (01.07.2012 -- 30.06.2013).}
\label{fig:demand}
\end{figure}

Furthermore, we set the observed period of one year from July 2012 to July 2013 and used the records of the German Meteorological Service for the weather station Bremen (airport BRE) as the data basis for the climate conditions to assign each day with a certain type of weighting factor, contrary to what is envisaged in the VDI~4655. The annual heat demand for the chosen household is 10741.15~kWh/a; this is about 6~$\%$ less than the theoretical guideline value. For a higher resolution, the total demand values for each day are further divided equally to observe the hourly demands depending on the daily energy requirements. From the preliminary considerations we calculated an example load profile for the following household: single-family household (short:~SFH) with 2 persons, living space of 100~m$^2$, specific heat demand of 100~kWh/(m$^2 \cdot $a) and an annual heat demand of 10741.15~kWh/a, which is illustrated in figure~\ref{fig:demand}.

\begin{figure}
\centering
\includegraphics[width=\textwidth]{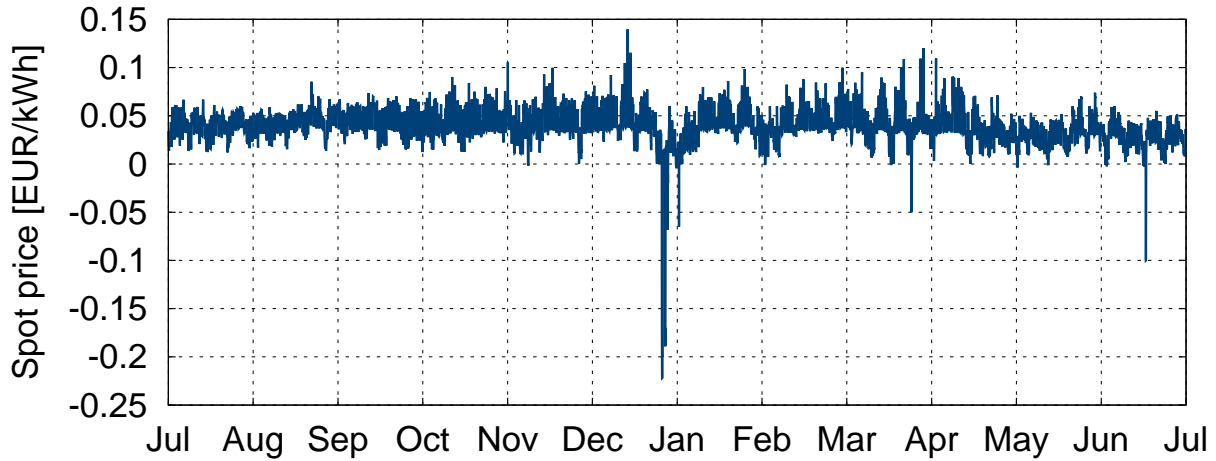}
\caption{Hourly price development of the electricity price (PHELIX) at the spot market EPEX SPOT, (01.07.2012 -- 30.06.2013).}
\label{fig:phelix}
\end{figure}

As a data basis for the electricity acquisition costs representing the volatility of the RES feed-in, we used the Physical Electricity Index (PHELIX). This reference price signal for the observed period published on the website \emph{http://www.epexspot.com} for auction results of the PHELIX and is illustrated in figure~\ref{fig:phelix} with an average price value of about 0.04~EUR/kWh.

\subsection{Cost-optimal charging strategies}
\label{subsec:subsec1}

We consider the SFH for two different PtH storage systems in the following. To evaluate and compare the examples more efficiently and even better, we set the maximal storage capacities for the thermal storage units as relative numbers corresponding to the mean daily average heat demand. The different observed systems are first a common PtH storage system consisting of a half-day storage unit (\(\approx\)~14.71~kWh) and a PtH conversion unit with a maximal charge power of 9~kW and second a system defined by a 14-day storage (\(\approx\)~411.99~kWh) and a conversion unit with a maximal charge power of 25~kW. This system neglects the switching limits of today's power connections and the today's common dimension of heat storage units in private households, but is chosen to demonstrate the cost saving potentials and also the overall benefit of a cost-optimal operation in combination with systems beyond today's standards. As mentioned above, the efficiencies of both conversion units for the moment are set to 100~$\%$. These observed systems are denoted as SFH$_{1/2}$ and SFH$_{14}$ respectively. As a reference value for the annual electricity acquisition costs and in order to compare the different cost savings, the costs with no storage unit would be 443.13~EUR to cover the annual heat demand; this represents the conversion of electricity into heat at times when the energy is required.

\begin{figure}
\centering
\includegraphics[width=\textwidth]{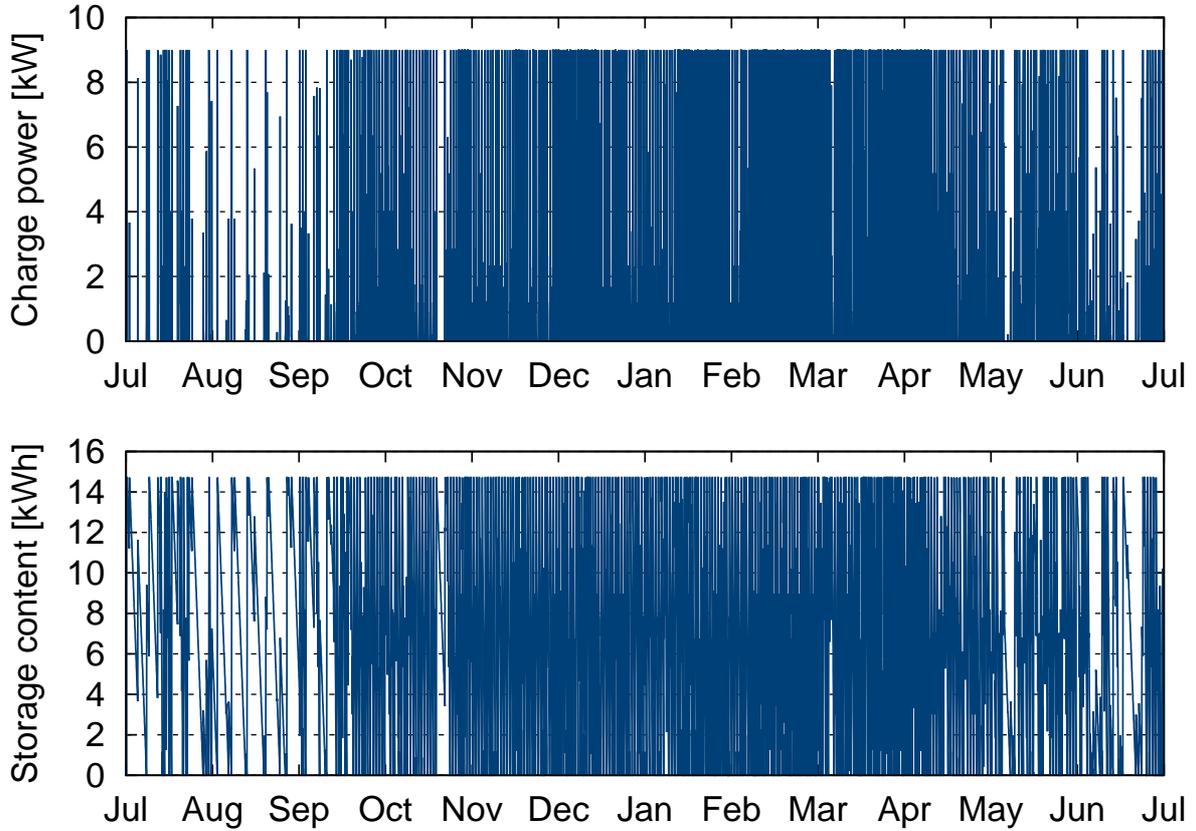}
\caption{Cost-optimal hourly charge curve (top panel) and corresponding hourly storage level (bottom panel) for a storage unit with capacity 14.71~kWh, (01.07.2012 -- 30.06.2013). The solution is 310.99~EUR for the electricity acquisition costs. (Assumptions: single-family household with 2 residents, living area of 100~m\(^2\), specific heat demand of 100~kWh/(m$^2 \cdot $a), storage capacity of 14.71~kWh (1/2-day of average heat demand) and charge power of 9~kW, no thermal energy losses)}
\label{fig:example1day}
\end{figure}

\begin{figure}
\includegraphics[width=\textwidth]{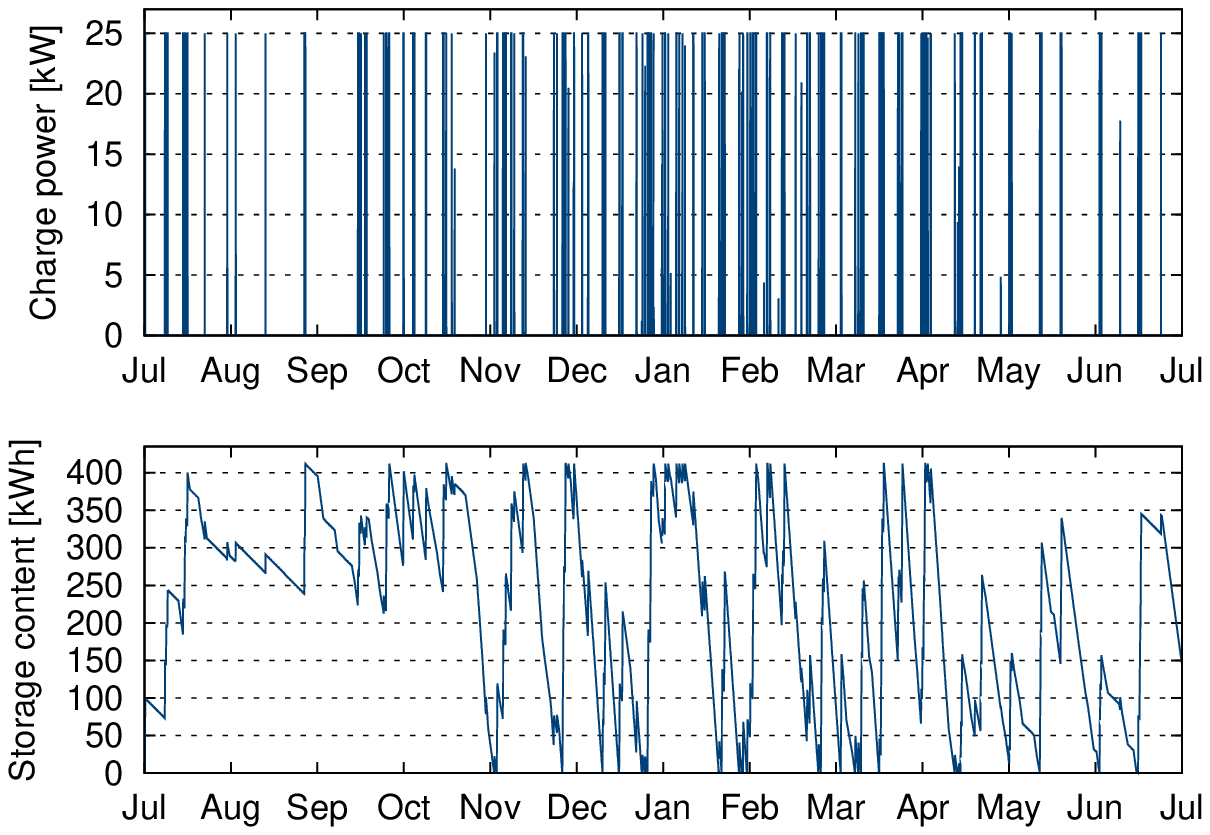}
\caption{Cost-optimal hourly charge curve (top panel) and corresponding hourly storage level (bottom panel) for a storage unit with capacity 411.99~kWh, (01.07.2012 -- 30.06.2013). The solution is 71.55~EUR for the electricity acquisition costs. (Assumptions: single-family household with 2 residents, living area of 100~m\(^2\), specific heat demand of 100~kWh/(m$^2 \cdot $a), storage capacity of 411.99~kWh (14-day of average heat demand) and charge power of 25~kW, no thermal energy losses)}
\label{fig:example14day}
\end{figure}

The calculated curves in figure~\ref{fig:example1day} for SFH$_{1/2}$ and in figure~\ref{fig:example14day} for SFH$_{14}$ result from the new algorithm and describe the cost-optimal strategy for operating the modeled PtH storage unit without energy losses. For SFH$_{1/2}$, the overall annual electricity acquisition costs are 310.99~EUR, which is 29.82~$\%$ less than the annual costs resulting from a system with no storage unit. The overall costs for SFH$_{14}$ are 71.55~EUR, which is more than 83~$\%$ below the reference value. The increased heat demand between November and April (cf. figure~\ref{fig:demand}) is represented by a higher number of charges within this period to cover the demand and by a highly fluctuating storage level. This behaviour is much more pronounced in the result of SFH$_{1/2}$ compared with SFH$_{14}$, which is due to the smaller storage capacity and also the lower charge power (cf. figures \ref{fig:example1day} and \ref{fig:example14day}). As described, the storage unit is charged as much as possible during hours with negative prices. This leads to a storage level greater or equal to zero at the end of the observed period. The storage level of SFH$_{1/2}$ is zero at the end of the period (cf. figure~\ref{fig:example1day}), compared to a storage level of 146.32~kWh for SFH$_{14}$ (cf. figure \ref{fig:example14day}). As a consequence, the large PtH storage system utilizes the negative electricity prices considerably better than the common smaller system.

\begin{figure}
\centering
    \begin{subfigure}[b]{.5\textwidth}
       \centering
       \includegraphics[width=\textwidth]{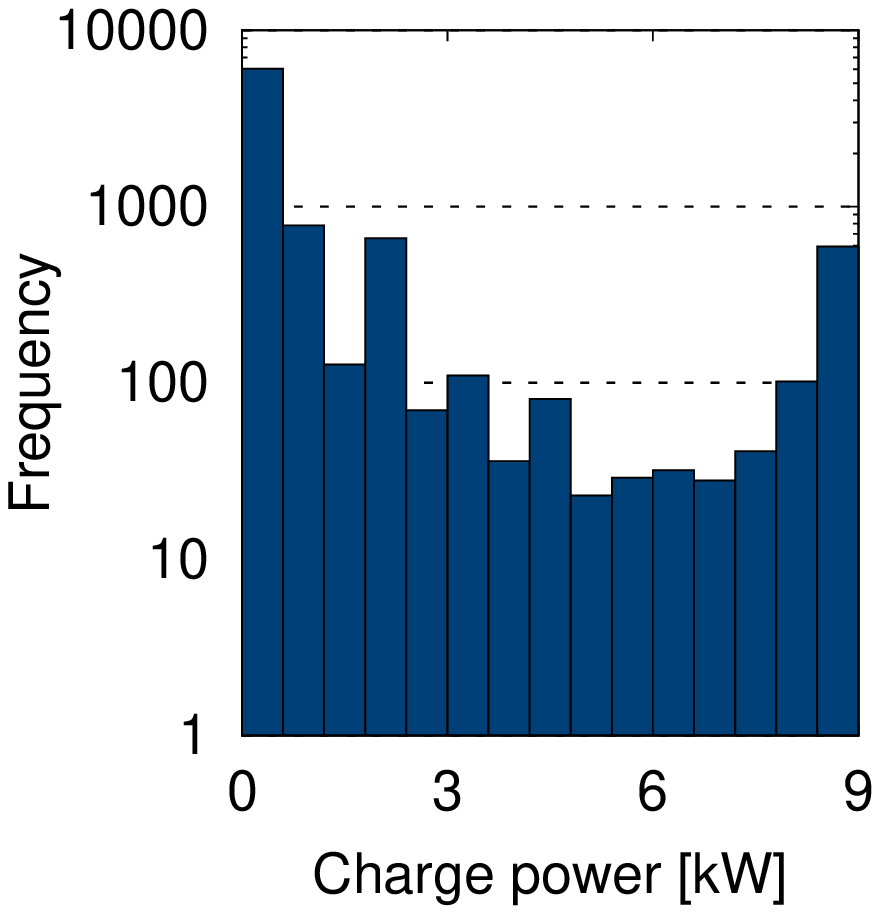}
       \caption{}
       \label{fig:hist1}
    \end{subfigure}
    \hspace{-5mm}
    \begin{subfigure}[b]{.5\textwidth}
       \centering
       \includegraphics[width=\textwidth]{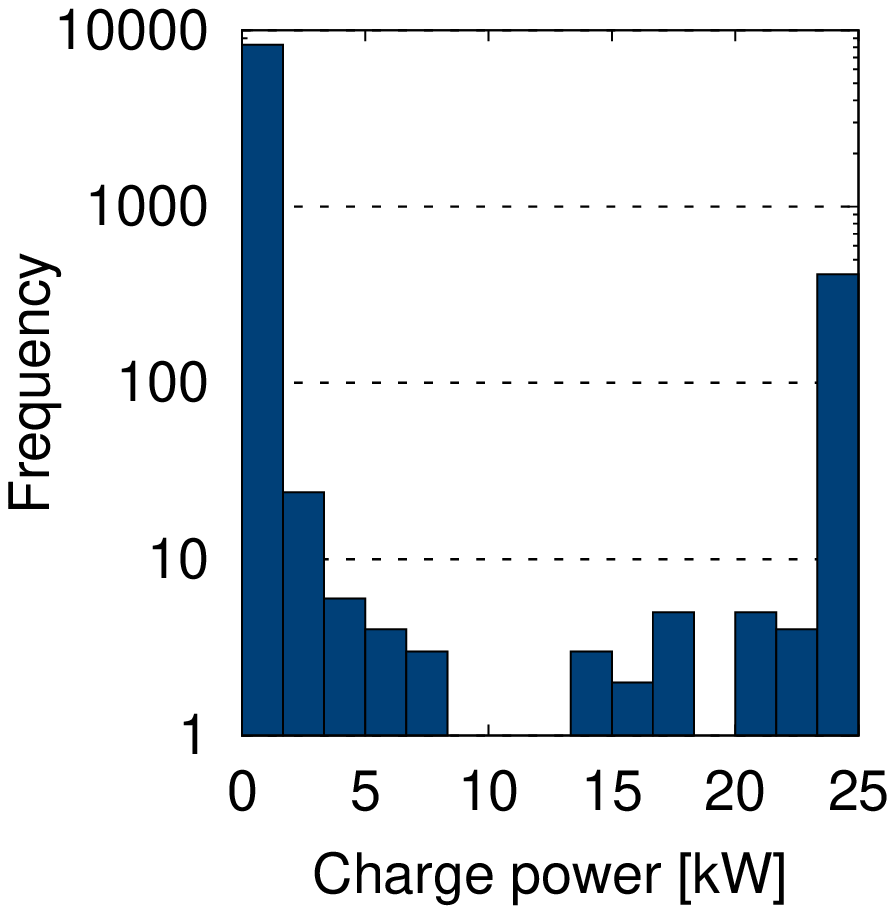}
       \caption{}
       \label{fig:hist14}
   \end{subfigure}
   \caption{(a) Histogram of the cost-optimal charges with a maximum charge power of 9~kW and storage capacity of 14.71~kWh (1/2-day) and an efficiency of 100~$\%$. (b) Histogram of the cost-optimal charges with a maximum charge power of 25~kW and storage capacity of 411.99~kWh (14-day) and an efficiency of 100~$\%$.}
   \label{fig:histogram}
\end{figure}

The difference between both storage capacities is also reflected in the decreasing number of non-zero charges from 2872 for SFH$_{1/2}$ down to 514 for SFH$_{14}$ as the capacity and charge power increases. For SFH$_{1/2}$, the total number of charges with maximum charge power is 575 and more in comparison to 407 maximum charges for SFH$_{14}$ (cf. figure~\ref{fig:histogram}). The number of charges, which are non-zero and less than the maximum charge power, turn out to be very different (cf. figure~\ref{fig:histogram}). These are 2297 for SFH$_{1/2}$ compared with 107 (less than 5~$\%$) for SFH$_{14}$. In summary, the storage capacity and the charge power have a strong impact on the behaviour of the resulting cost-optimal charge curve for the chosen household. As these factors rises, the number of charges lower than the maximum charge power decreases drastically and the number of charges with the maximum power increases slightly. This leads to a shift in the values of the cost-optimal charge strategy towards the maximum charge power. This fact is presented in figure~\ref{fig:hist1} for SFH$_{1/2}$ and in figure~\ref{fig:hist14} for SFH$_{14}$. Further, increasing the system factors reduces operating hours tremendously, which leads in combination with a storage level higher than zero at the end to a sharply higher cost saving potential for large systems.

We now add energy losses of the thermal storage unit to the problem as described in section~\ref{sec:sec2} and further maintain no energy losses of the conversion system. With respect to the directive 2010/30/EU of the European Parliament and of the Council \cite{directive2010}, we applied an exponential regression analysis on the data taken from \cite{vaillant2015}, to calculate heat losses in standby mode (details in~\ref{app:regression}). The result of the regression leads to an hourly factor of about 0.9962 (approx. 8.73~$\%$ daily energy losses) for the half-day storage unit and an hourly factor of about 0.9989 (approx. 2.61~$\%$ daily energy losses) for the 14-day storage unit.

\begin{figure}
\includegraphics[width=\textwidth]{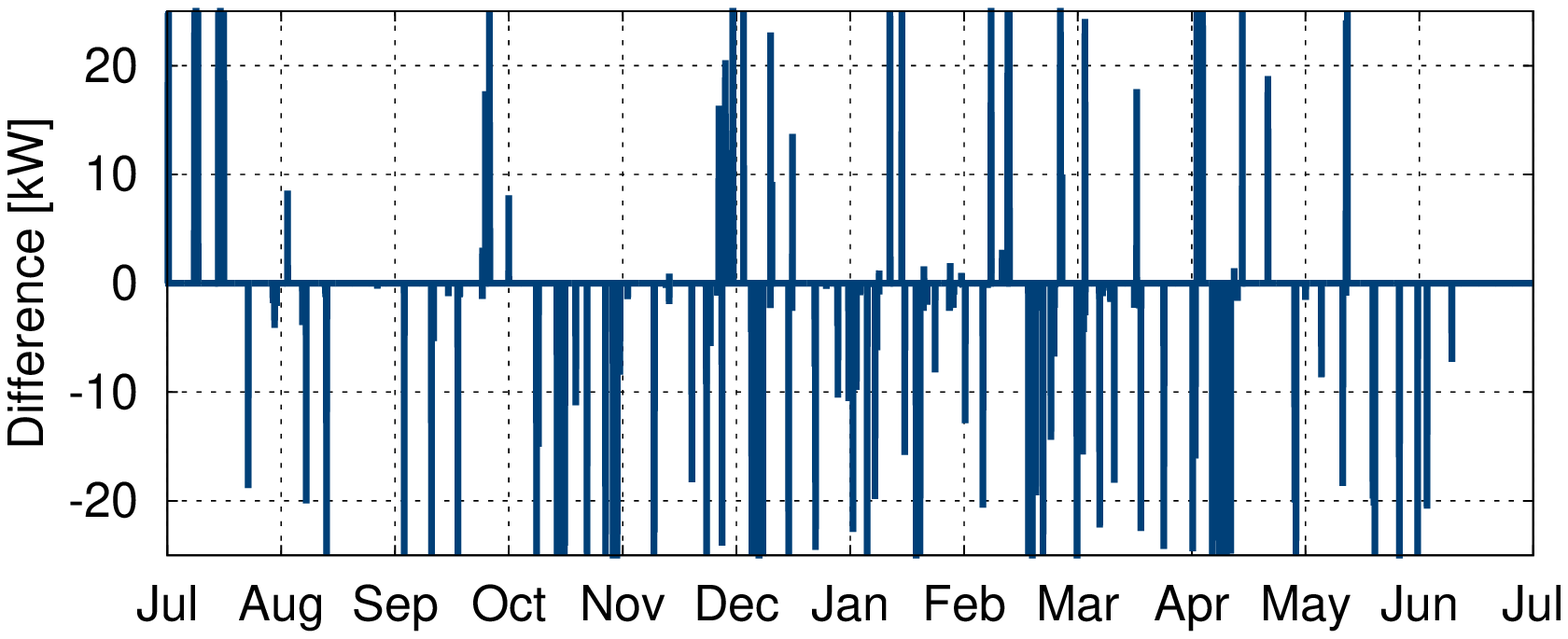}
\caption{Differences between both cost-optimal charging strategies for a storage unit with a capacity of 411.99~kWh (14-day) and hourly efficiencies of 100~$\%$ and 99.89~$\%$, respectively. Positive values indicate higher charges at the respective times in the 100~$\%$-case. As a general trend decreased efficiencies result in increased energy consumption and delayed charging of the storage. The additional annual amount of energy required for cost-optimal charging of the storage unit with finite efficiency of 99.89~$\%$ is 994.19~kW.}
\label{fig:diff_14day}
\end{figure}

The differences between the cost-optimal charge curve without energy losses and the solution of the problem with losses for SFH$_{1/2}$ are quite small. The charges of the different charge curves are shifted by only a few hours in time, which means that the moment in time and the amount of energy are almost identical for both examples. As a result, almost no differences between the cost-optimal charge curves of a heat storage system with a lossless and a lossy small storage capacity. By contrast, the differences in the SFH$_{14}$ versions display a block-like structure. The behaviour is either maximum charge power for the cost-optimal charges with a lossless storage unit or for the lossy version over a couple of weeks (cf. figure \ref{fig:diff_14day}). The annual electricity acquisition costs for SFH$_{1/2}$ with a lossy storage unit is 318.09~EUR and 96.78~EUR for the lossy version SFH$_{14}$. This means an increase of about 2.28~$\%$ for the first and 35.26~$\%$ for the second system. Nevertheless, it is 28.22~$\%$ and 78.16~$\%$ respectively less than the reference value. This is caused by the additional amount of energy needed to cover the heat demand, which is 200.78~kWh for SFH$_{1/2}$ and 994.19~kWh for SFH$_{14}$ (cf. figure \ref{fig:diff_14day}).

\begin{figure}
\centering
\includegraphics[width=\textwidth]{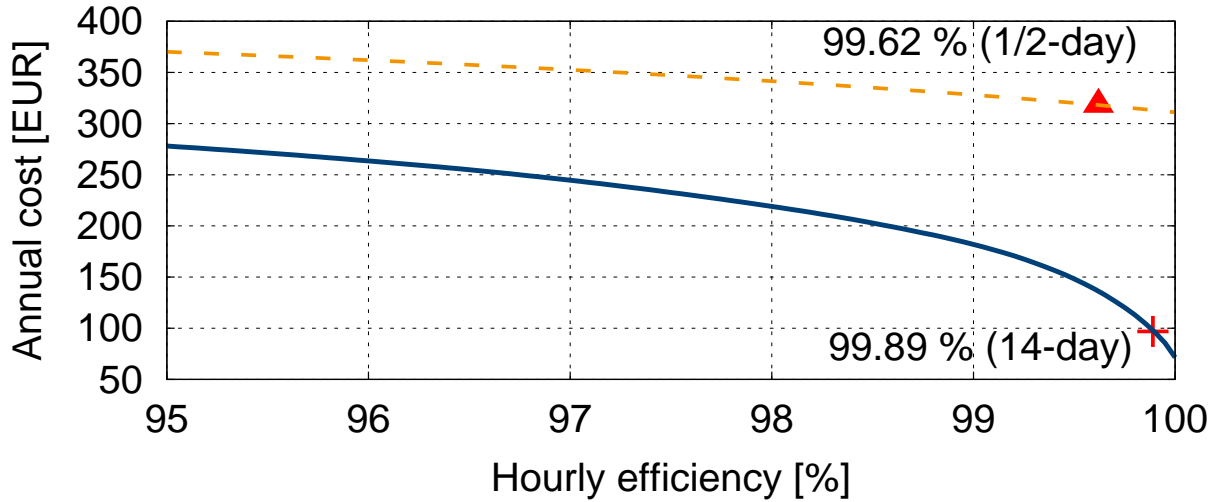}
\caption{Electricity acquisition costs for PtH storage systems with a maximum charge power of 9~kW and a storage capacity of half a day ($\approx$~14.71~kWh, dashed orange line) and also 25~kW charge power in combination with a storage capacity of 14 days ($\approx$~411.99~kWh, solid blue line), both for hourly storage efficiencies from 95~$\%$ up to 100~$\%$. The calculated factors of energy losses for these storage capacities are about 99.62~$\%$ (red triangle) for the 1/2-day storage unit, which leads to overall costs of 318.09~EUR and 99.89~$\%$ (red cross) for the 14-day storage unit, which results in overall costs of 96.78~EUR.}
\label{fig:efficiency factor}
\end{figure}

Below we address the question of how the efficiency factor affects the solution. For the two example systems with a charge power of 9~kW and 25~kW respectively and two different storage capacities of half a day and 14 days as well as the household specifications as above, the results of a variable hourly efficiency of the thermal storage unit for both examples are presented in figure~\ref{fig:efficiency factor}. The hourly efficiencies are between 0.95 and 1.00, which represent daily energy losses of between 70~$\%$ and 0~$\%$. The annual electricity acquisition costs decrease as the efficiency increases in both cases. This is due to the lower amount of energy required for an improved isolated storage unit. The curvature of the curve for the 14-day storage unit is growing faster as the efficiency increases, in comparison with the curve for the half-day storage unit. This leads to higher influences of the efficiencies for larger storage units, which have been indicated in the previous examples. In the case of hourly efficiencies above 0.99 (daily energy losses of about 22~$\%$) in particular (cf. figure \ref{fig:efficiency factor}), the optimal annual costs for larger capacities decrease more quickly compared with the trajectory of small ones.

\subsection{Operating costs for different system combinations}
\label{subsec:subsec2}

Now we present the application of sensitivity analysis to the optimization problem and point out the benefits of the new algorithm. The considerations in \ref{subsec:subsec1} were based on given (fixed) values C and S for the technical restrictions on the conversion unit and storage capacity. In addition to the acquisition costs, installation costs for a PtH storage unit will be observed from now on. The aim is to minimize the annual electricity acquisition costs, with the inclusion of costs for the purchase and installation of different PtH storage units. This more economic approach for estimations of the overall costs leads to the minimization problem:

\begin{equation}\label{eq:nonlinear problem}
\min_{\mathcal{C},\mathcal{S} \le \infty} \mathcal{K}\left(\mathcal{C}\right) + \mathcal{L}\left(\mathcal{S}\right) + \mathcal{M}\left(\mathcal{C},\mathcal{S} \right)\,,
\end{equation}

\noindent with (non-linear) functions \(\mathcal{K}\left(\mathcal{C}\right)\), \(\mathcal{L}\left(\mathcal{S}\right)\) and \(\mathcal{M}\left(\mathcal{C},\mathcal{S} \right)\) as the initial minimization problem~\eqref{eq:LP}. Here, \(\mathcal{K}\left(\mathcal{C}\right)\) represents the cost-function of a PtH unit depending on a certain maximum charge power $\mathcal{C}$ and \(\mathcal{L}\left(\mathcal{S}\right)\) the cost-function for a storage unit with capacity $\mathcal{S}$. The solution of \eqref{eq:nonlinear problem} requires $\mathcal{C}$ times $\mathcal{S}$ calculations of the optimization problem \eqref{eq:P}. Due to the gain in efficiency of the new algorithm (see:~\ref{app:algorithm}), solving the problem \eqref{eq:LP} for even thousands of different combinations \(\mathcal{C}, \mathcal{S}\) is possible within short time frame (cf. table~\ref{table:runtime}).

\begin{figure}
\centering
\includegraphics[width=.8\textwidth]{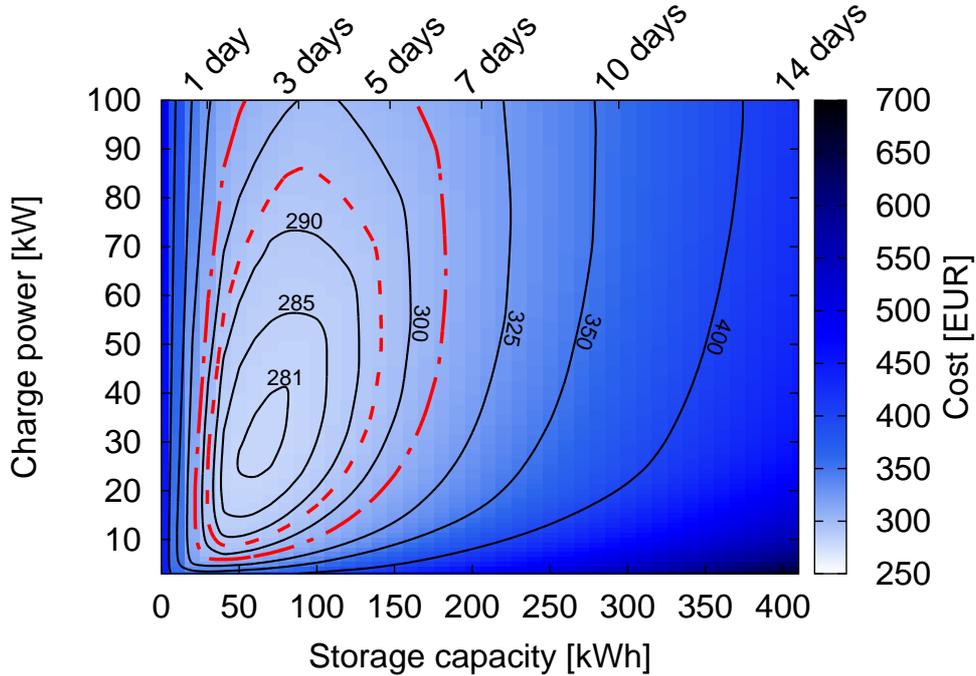}
\caption{Total costs for different system setups. The minimum value is 280.02~EUR for a system with \(\mathcal{C}=28\)~kW and \(\mathcal{S}=70\)~kWh (about two and a half days). The 5~$\%$ level-line (red dashed line) is at 294.02~EUR and the 10~$\%$ level-line (red solid-dashed line) at 308.02~EUR. (Assumptions: life span of 15 years, charge power between 3~kW and 100~kW, storage capacity between 0~kWh (no storage unit) and 410~kWh (14-day storage unit) and an efficiency of 100~$\%$)}
\label{fig:contour1}
\end{figure}

The household data set, the observed time-period and the price signal are the same as in \ref{subsec:subsec1} and the life span of the installed PtH storage system is assumed to be 15 years. For the time being, the efficiencies of the storage units are set to 100~$\%$. For ease of presentation maintenance costs and sales taxes are not included in the analysis at this stage.

Based on regressions detailed in~\ref{app:regression}, the partial cost-functions \(\mathcal{K}\left(\mathcal{C}\right)\) and  \(\mathcal{L}\left(\mathcal{S}\right)\) for the annual costs of charging and storage units can be approximated from data, resulting in the specific  optimization problem
\begin{equation}\label{eq:nonlinear problem example1}
\min \ \left(0.47 \cdot \mathcal{C} \ \text{ EUR/kW}\right) + \left(0.95 \cdot \mathcal{S} \ \text{ EUR/kWh}\right) + \mathcal{M}\left(\mathcal{C},\mathcal{S} \right)\,.
\end{equation}
The maximum hourly demand of the example household is about 3~kWh, so we set \linebreak \(3 \le \mathcal{C} \le 100\), where the upper bound may not be feasible in practice and only represents a theoretical value to analyze the behaviour of the solutions. The observed storage capacities are between zero, which represents a system without any storage unit, and a storage capacity of 410~kWh, which is about the 14-day daily mean average heat demand of the household. The subdivision of \(\mathcal{S}\) is set to 10~kWh (\(\approx\) 8~hours) and the separation of \(\mathcal{C}\) is set to 1~kW.

The contours of the solutions for each combination are illustrated in figure~\ref{fig:contour1}. The solution of the minimization problem~\eqref{eq:nonlinear problem example1} is 280.02~EUR for the combination \(\mathcal{C}=28\)~kW and \(\mathcal{S}=70\)~kWh (cf. figure \ref{fig:contour1}). It has to be mentioned that this solution does not represent additional installation, operational and further fix costs as discussed above. In terms of sensitivity of the solution to the problem~\eqref{eq:nonlinear problem example1}, we study the area surrounding the minimum. The 5~$\%$ and 10$\%$ levels, which represent the area of solutions that are 5~$\%$ and 10$\%$ respectively higher than the optimum, are 294.02~EUR and 308.02~EUR (cf. fig.~\ref{fig:contour1}). This comprises a wide range of the calculated points and is about 15.74~$\%$ of all calculated points for the 5~$\%$ level as well as more than 30~$\%$ for the 10~$\%$ level. These areas, which represent about the same overall costs, can be used to estimate economic combinations for PtH systems and connected storage units. A result of this approach is that the combination for a charge power of 6~kW and a storage capacity of 30~kWh results in about the same total costs for one year as a combination of 100~kW for the charge power and 160~kWh for the storage capacity. This leads to a relatively robust and insensitive behaviour of the solution, but depends on the functions \(\mathcal{K}\left(\mathcal{C}\right)\) and \(\mathcal{L}\left( \mathcal{S}\right)\) in~\eqref{eq:nonlinear problem example1}. Further, figure~\ref{fig:contour1} shows that the inclusion of installation costs reduces the cost savings potentials especially for large storage systems, in contrast to the results in~\ref{subsec:subsec1}.

\subsection{Operating costs for different system combinations including losses}
\label{subsec:subsec3}

\begin{figure}[t]
\centering
\includegraphics[width=.8\textwidth]{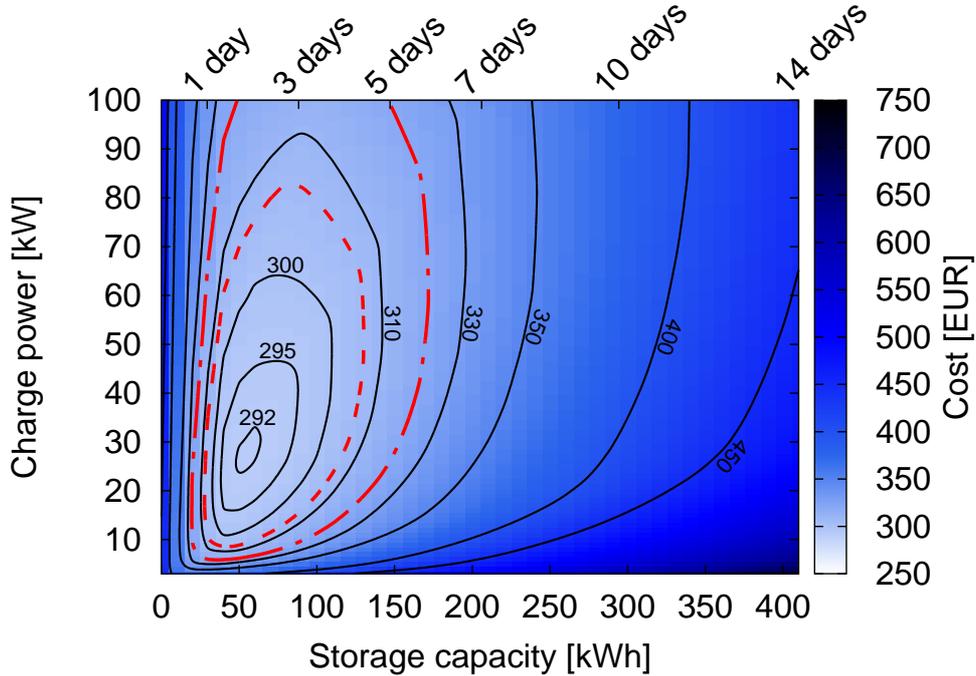}
\caption{Total costs for different systems with the same representation and assumptions as described already in figure \ref{fig:contour1}, but now with thermal storage losses explicitly included in the analysis (see section \ref{subsec:subsec3} for details). The resulting cost-minimum is 291.50~EUR for a system with \(\mathcal{C}=26\)~kW and \(\mathcal{S}=60\)~kWh (about two days). The 5~$\%$ level-line (red dashed line) is at 306.08~EUR and the 10~$\%$ level-line (red solid-dashed line) at 320.65~EUR.}
\label{fig:efficiency contour}
\end{figure}

In~\ref{subsec:subsec2} we considered the storage units with an efficiency of 100~$\%$. We now add energy losses of the storage unit to the problem \eqref{eq:nonlinear problem example1}, which will be dependent on the maximum storage capacity. Therefore, the function \(\mathcal{M}\left( \mathcal{C},\mathcal{S} \right)\) in problem~\eqref{eq:nonlinear problem example1} is now represented by~\eqref{eq:LP variable}. We use the same procedure as described in~\ref{app:regression} to calculate hourly energy losses for each storage capacity. The data sets are assumed to be the same as in \ref{subsec:subsec2}, but with additional loss factors for the different storage capacities $\mathcal{S}\in \{0~\text{kWh},10~\text{kWh},20~\text{kWh}, \dots ,410~\text{kWh}\}$. The contours of the solutions are illustrated in figure \ref{fig:efficiency contour}. 

The solution of the minimization problem with lossy storage units is 291.50~EUR for the combination \(\mathcal{C}=26\)~kW and \(\mathcal{S}=60\)~kWh (cf. figure \ref{fig:efficiency contour}), which is just 4.1~$\%$ higher than the solution in~\ref{subsec:subsec2}. Again, we study the sensitivity and robustness of the solution. The 5~$\%$ level has a value of about 306.08~EUR and the 10~$\%$ level is 320.65~EUR. As a result, the consideration of storage units with energy losses leads to a solution which is close to the solution of \eqref{eq:nonlinear problem example1}, but the area with equal solutions for different combinations is smaller. Therefore, the overall costs for combinations with energy losses increase faster away from the minimum than the costs without losses, but are still relatively insensitive.

\section{Conclusions}
\label{sec:sec7}

The integration of high shares of electricity generated by RES into existing electricity grids constitutes a major challenge for future energy systems. Energy storage units can contribute significantly to the expansion and integration of electricity from RES.

In section~\ref{sec:sec2}, we developed a mathematical model to describe conversion systems combined with storage units. Based on the idea of using electricity from RES efficiently, utilizing local overproduction and increasing the share of RES within electricity grids, we assumed that the electricity price represents the availability of electrical energy as well as the volatile feed-in from RES. Next, we investigated the cost-optimal charging of these energy storage units, which is linked to the minimum electricity acquisition costs. This minimization problem was expressed as a mathematical optimization problem. Furthermore, a procedure for the inclusion of constant and non-constant energy losses was described, which covers a wide range of applications.

The focus of this work was the development of a new optimization method to calculate solutions for the modeled storage-related optimization problems and proof of optimality. Due to the special structure of the storage problem, this method provides a tremendous gain in efficiency and is also applicable to general storage-related optimization problems.

In section \ref{sec:sec3}, we applied the approach from the previous section to the PtH technology combined with thermal storage units. The aim was to study the principal effects of our approach and the consequences on the resulting cost-optimal charging strategies. Based on a fluctuating electricity price as a measure for the electricity acquisition costs, we calculated cost-optimal charging strategies for different PtH storage systems including storage units with and without energy losses. We found that the increase in capacity leads to a clustering of the cost-optimal charge values. Most of them are either zero or maximum charge power, which leads to fewer fluctuations of the storage level over the observed time period. Furthermore, we evaluated the impact of a variable efficiency factor on the minimum electricity acquisition costs for two different storage capacities, resulting in a higher efficiency impact for large storage units. We also described an approach to estimate cost-optimal combinations of the conversion system and the storage unit (with and without energy losses), with the inclusion of installation costs. The result is that the inclusion of energy losses decreases the cost-optimal solution only by a small amount and retains the elliptic shape compared with the lossless solutions. This leads to limited influences of different efficiencies on the cost-optimal combination for the observed household. But the inclusion of installation costs significantly decreases the cost savings potentials especially for large PtH storage systems. At the end of this section, we illustrated the robustness of the overall costs with respect to the optimal solution, as a consequence of the applied sensitivity analysis and have shown that this optimal solution can be reached by systems, which are larger than today.

In this work, we used the fluctuating electricity spot market price as a measure for the approximation of the electricity acquisition costs. In practice, this measure underestimates common acquisition costs. On the contrary, the utilization of a conversion process with a high performance factor, e.g. a heat pump system, could have a compensatory effect. Nevertheless, the investigations in this work illustrate the resultant effects of a cost-optimal charging strategy.

In conclusion, the novel optimization method developed in this work is applicable to storage-related optimization problems, even where different efficiencies are included. Due to the reduction in runtime, the presented method can also be used to manage a collection of decentralized storage units, for instance via the utilization of small ECUs. Furthermore, this work reveals the benefits of a problem-specific approach in terms of mathematical simplicity, applicability to general storage-related problems and decrease of computational effort for even high-resolution problems on a large time-scale.

\section{Acknowledgements}
The authors kindly acknowledge preliminary research and a contribution to the preparation of data on the heat demand of single-family houses by Konrad Meyer. This data was used for the application of the algorithm to PtH storage units presented in section~\ref{sec:sec3}. Funding of the joint project RESTORE 2050 (funding code 03SF0439A) was kindly provided by the German Federal Ministry of Education and Research through the funding initiative Energy Storage.


\appendix

\def \NN {\mathbb{N}}
\def \RR {\mathbb{R}}
\def \CC {\mathbb{C}}
\def \ZZ {\mathbb{Z}}
\def \KK {\mathbb{K}}
\def \QQ {\mathbb{Q}}
\def \Rcal {\mathcal{R}}
\def \Ncal {\mathcal{N}}
\def \Lcal {\mathcal{L}}
\def \Ocal {\mathcal{O}}
\newcommand{\set}[2]{\{#1\, |\, #2\}}

\renewcommand{\algorithmicrequire}{\textbf{Input:}}
\renewcommand{\algorithmicensure}{\textbf{Output:}}

\section{New algorithm: definition and proof of optimality}\label{app:algorithm}

Let $a_i,b_i,c_i,u_i \in \RR$ with $u_i \ge 0$ and $a_i \le b_i$ for $i\in I:=\{1,\ldots,n\}$.
We are interested in finding a point $x=(x_1,\ldots,x_n) \in \RR^n$ which solves the following constrained minimization problem~\eqref{eq:P},
\begin{align} 
\min \quad & \sum_{j \in I} c_j \cdot x_j   \nonumber \\
\mbox{subject to} \quad & 0 \le x_i \le u_i \label{eq:P} \tag{P} \\
\mbox{and}\quad  & a_i \le \sum_{j=1}^i x_j \le b_i \quad \mbox{for all $i\in I$} \,.  \nonumber
\end{align}
We denote its minimum value by $P_{min}$. A \emph{feasible point} of~\eqref{eq:P} is any $x \in \RR^n$ which fulfills the constraints, and a \emph{solution} of~\eqref{eq:P} is a feasible point $\hat{x}$ which achieves the minimum value $P_{\min}$, i.e. $\sum_{j \in I} c_j \cdot \hat{x}_j = P_{\min}$. Note that, although $P_{\min}$ is unique, there could be different solutions $\hat{x}$ if some of the $c_i$ coincide or are equal to zero. We assume that there exists at least one feasible point. Hence~\eqref{eq:P} has at least one solution, because the constraints $0 \le x_i \le u_i $ ensure that the set of all feasible points is compact. For more information about the solvability of linear programs we refer to~\cite{nocedal2006numerical}. There the reader can also find common algorithms to solve linear programs, where among the best known are the \emph{Simplex Method} and \emph{Interior Point Methods}. Here we exploit the special structure of~\eqref{eq:P} to develop the much more efficient new algorithm, which needs at most $n^2+3n$ floating point operations and \(3/2n^2+5/2n\) comparisons to compute a solution. The source code of this algorithm written in the programming language Python is presented in \ref{app:python}.

\begin{algorithm}[ht]
\renewcommand{\thealgorithm}{}
  \caption{}\label{alg}
    \begin{algorithmic}[1]
        \Require{$a_i,b_i,c_i,u_i$ and the permutation $\sigma$}
        \Ensure{A solution \(x\) of the problem~\eqref{eq:P}}
        \For{\(k=1\) to \(n\)}
	        \State \(M_1 \gets \max\limits_{i<\sigma(k)}\{0,a_i\}\)
            \State \(M_2 \gets \max\limits_{i \ge \sigma(k)}\{0,a_i\}\)
            \State \(m \gets \min\limits_{i \ge \sigma(k)} \{b_i\}\)
	        \If{\(c_{\sigma(k)} \ge 0\)}
				\State $x_{\sigma(k)} \gets \min\big\{\max\{0,M_2-M_1\},\min\{u_{\sigma(k)}, m-M_1\}\big\}$
			\Else
				\State $x_{\sigma(k)} \gets \min\{u_{\sigma(k)}, m-M_1\}$
			\EndIf
			\For{$i=\sigma(k),\ldots,n$}
		    \State $a_i \gets a_i-x_{\sigma(k)}$
		    \State $b_i \gets b_i-x_{\sigma(k)}$
	    \EndFor
    \EndFor
  \end{algorithmic}
\end{algorithm}

We now give a precise description of the new algorithm, which is denoted as Algorithm in the following. Let $\sigma$ be a permutation of the indices which correspond to increasing $c_i$ values, i.e. $c_{\sigma(1)} \le \ldots \le c_{\sigma(n)}$, and let $\hat{x}$ be a solution of~\eqref{eq:P}. Then Algorithm successively computes the values $\hat{x}_{\sigma(1)}, \ldots, \hat{x}_{\sigma(n)}$. In fact we show in the following 5 steps that the value $\hat{x}_{\sigma(k)}$ must necessarily be equal to the explicitly computable value in Line 6 or 8 of Algorithm, which proves that Algorithm indeed computes a solution of~\eqref{eq:P}. At first we assume that all $c_i$ are different, i.e. $c_{\sigma(1)} < \ldots < c_{\sigma(n)}$, and not equal to zero.

\noindent \textbf{Step 1}: We show that $P_{min}=c_{\sigma(1)} \cdot\hat{x}_{\sigma(1)} + P_{min}'$, where $P_{min}'$ is the minimum value of the reduced linear program
\begin{align} 
\min \quad &\sum_{j=2}^n c_{\sigma(j)} \cdot x_{\sigma(j)}    \nonumber \\
\mbox{subject to} \quad  & 0 \le x_i \le u_i \label{eq:Pprime} \tag{P'} \\
\mbox{and}\quad  & a_i' \le \sum_{\genfrac{}{}{0pt}{}{j=1}{j \neq \sigma(1)}}^i x_j \le b_i' \quad \mbox{for all $i \in I \setminus \{\sigma(1)\}$}\,,  \nonumber
\end{align}
with the modified lower and upper bounds
\begin{align*} 
a_i' &:= \begin{cases} a_i &, i < \sigma(1)-1 \\ \max\{a_{\sigma(1)-1},a_{\sigma(1)}-\hat{x}_{\sigma(1)}\} &, i = \sigma(1)-1 \\ a_i- \hat{x}_{\sigma(1)} &, i>\sigma(1) \end{cases} \\
b_i' &:= \begin{cases} b_i &, i < \sigma(1)-1 \\  \min\{b_{\sigma(1)-1},b_{\sigma(1)}-\hat{x}_{\sigma(1)}\} &, i = \sigma(1)-1 \\b_i- \hat{x}_{\sigma(1)} &, i>\sigma(1) \end{cases}
\end{align*}
Indeed, the point $(\hat{x}_1,\ldots,\hat{x}_{\sigma(1)-1},\hat{x}_{\sigma(1)+1},\ldots,\hat{x}_n) \in \RR^{n-1}$ is obviously feasible for~\eqref{eq:Pprime}, and hence we have
\[
P_{min}=c_{\sigma(1)} \cdot\hat{x}_{\sigma(1)} + \sum_{j=2}^n c_{\sigma(j)} \cdot \hat{x}_{\sigma(j)} \ge c_{\sigma(1)} \cdot\hat{x}_{\sigma(1)} + P_{min}'\,.
\]
Conversely let $(x_1,\ldots,x_{\sigma(1)-1},x_{\sigma(1)+1},\ldots,x_n) \in \RR^{n-1}$ be a solution of~\eqref{eq:Pprime}.
Then the point $(x_1,\ldots,x_{\sigma(1)-1},\hat{x}_{\sigma(1)},x_{\sigma(1)+1},\ldots,x_n) \in \RR^n$ is feasible for~\eqref{eq:P}, and hence we also have
\[
c_{\sigma(1)} \cdot\hat{x}_{\sigma(1)} + P_{min}'=c_{\sigma(1)} \cdot\hat{x}_{\sigma(1)} + \sum_{j=2}^n c_{\sigma(j)} \cdot x_{\sigma(j)} \ge P_{min}\,.
\]
From this we moreover infer that $(\hat{x}_1,\ldots,\hat{x}_{\sigma(1)-1},\hat{x}_{\sigma(1)+1},\ldots,\hat{x}_n)$ is actually a solution of~\eqref{eq:Pprime}.
As a consequence, once the value of $\hat{x}_{\sigma(1)}$ is known, it remains to compute the remaining values $\hat{x}_{\sigma(2)},\ldots,\hat{x}_{\sigma(n)}$ by solving~\eqref{eq:Pprime}.
Since~\eqref{eq:Pprime} has the same structure as~\eqref{eq:P}, we can proceed analoguously to get $P_{min}'=c_{\sigma(2)} \cdot\hat{x}_{\sigma(2)} + P_{min}''$, with a similar problem (P''), and determine the value of $\hat{x}_{\sigma(2)}$.
Iterating this process $n$ times leads to the outer for-loop in Algorithm.
The inner for-loop corresponds to the modification of the lower and upper bounds in~\eqref{eq:Pprime} (note that only maxima of the $a_i$ and minima of the $b_i$ are needed to compute $\hat{x}_{\sigma(k)}$).
It remains to show that $\hat{x}_{\sigma(1)}$ equals the explicit value in Line 6 or 8 of Algorithm, because then this inductively holds for all $\hat{x}_{\sigma(k)}$.
As in Line 2--4 of Algorithm we define
\[
M_1:=\max_{i<\sigma(1)}\{0,a_i\} \quad, \quad M_2:=\max_{i \ge \sigma(1)}\{0,a_i\} \quad,\quad m:=\min_{i \ge \sigma(1)} \{b_i\} \,.
\]
\textbf{Step 2}: We show that $\sum_{j=1}^{\sigma(1)-1} \hat{x}_j \ge M_1$ and $0 \le \hat{x}_{\sigma(1)} \le \min\{u_{\sigma(1)},m-M_1\}$.
The first inequality holds because for all $i<\sigma(1)$ we have
\[
\sum_{j=1}^{\sigma(1)-1} \hat{x}_j \ge \sum_{j=1}^i \hat{x}_j \ge a_i \,.
\]
And the second inequality holds because of $0 \le \hat{x}_{\sigma(1)} \le u_{\sigma(1)}$, and because for all $i \ge \sigma(1)$ we have
\[
b_i \ge \sum_{j=1}^i \hat{x}_j \ge \hat{x}_{\sigma(1)} + \sum_{j=1}^{\sigma(1)-1} \hat{x}_j \ge \hat{x}_{\sigma(1)} +  M_1 \,.
\]
\textbf{Step 3}: We show that in case $\hat{x}_{\sigma(1)} < \min\{u_{\sigma(1)},m-M_1\}$ we must have $\sum_{j=1}^{\sigma(1)-1} \hat{x}_j=M_1$.
Assume to the contrary that $\sum_{j=1}^{\sigma(1)-1} \hat{x}_j>M_1 \ge 0$.
Then we find an index $i < \sigma(1)$ such that $\hat{x}_i>0$.
Let $i_1$ be maximal with this property, i.e. $\hat{x}_{i_1}>0$ and $\hat{x}_i=0$ for all $i_1 < i < \sigma(1)$.
Then we find some $\delta>0$ such that $\hat{x}_{i_1}- \delta>0$, $\sum_{j=1}^{\sigma(1)-1} \hat{x}_j - \delta > M_1$ and $\hat{x}_{\sigma(1)} +\delta < u_{\sigma(1)}$.
Now we define $x \in \RR^n$ by
\[
x_i:=\begin{cases} \hat{x}_{i_1}-\delta &, i=i_1 \\ \hat{x}_{\sigma(1)}+\delta &, i=\sigma(1) \\ \hat{x}_i &, \mbox{otherwise} \end{cases} \,.
\]
Then for all $i<i_1$ and all $i \ge \sigma(1)$ we have $\sum_{j =1}^i x_j = \sum_{j =1}^i \hat{x}_j$.
And for all $i_1 \le i < \sigma(1)$ we have 
\[
\sum_{j =1}^i x_j = \sum_{j =1}^i \hat{x}_j - \delta \quad \begin{cases} \quad \le b_i \\ \quad =\sum_{j =1}^{\sigma(1)-1} \hat{x}_j - \delta > M_1 \ge a_i \end{cases} \,.
\]
Hence $x$ is feasible for~\eqref{eq:P} with
\[
\sum_{j \in I} c_j \cdot x_j = P_{min} + \delta \cdot (c_{\sigma(1)} - c_{i_1}) < P_{min} \,,
\]
which is a contradiction to the fact that $\hat{x}$ is a solution of~\eqref{eq:P}.\\
\textbf{Step 4}: In case $c_{\sigma(1)}>0$ we prove that
\[
\hat{x}_{\sigma(1)}=\min\big\{\max\{0,M_2-M_1\},\min\{u_{\sigma(k)}, m-M_1\}\big\} \,.
\]
For this we consider three cases:
\begin{enumerate}
\item In case $\hat{x}_{\sigma(1)} = M_2-M_1$ it follows with the second inequality in Step 2 that $0 \le M_2-M_1 \le \min\{u_{\sigma(1)},m-M_1\}$, and hence the assertion is~true.
\item In case $\hat{x}_{\sigma(1)} < M_2-M_1$ we show that $\hat{x}_{\sigma(1)} = \min\{u_{\sigma(1)},m-M_1\}$, from which the assertion follows.
Assume to the contrary that $\hat{x}_{\sigma(1)} < \min\{u_{\sigma(1)},m-M_1\}$.
Then it follows from Step~3 that $\sum_{j=1}^{\sigma(1)-1} \hat{x}_j=M_1$.
Let $i_0 \ge \sigma(1)$ be an index with $a_{i_0}=\max_{i \ge \sigma(1)}\{a_i\}$.
Then we have
\[
M_2 \le \sum_{j =1}^{i_0} \hat{x}_j =  M_1 +\hat{x}_{\sigma(1)}+ \sum_{j =\sigma(1)+1}^{i_0} \hat{x}_j < M_2 + \sum_{j =\sigma(1)+1}^{i_0} \hat{x}_j\,,
\]
from which we infer that $\sum_{j =\sigma(1)+1}^{i_0} \hat{x}_j > 0$.
Hence we find an index $i > \sigma(1)$ such that $\hat{x}_i>0$.
Let $i_1$ be minimal with this property, i.e. $\hat{x}_{i_1}>0$ and $\hat{x}_i=0$ for all $\sigma(1) < i < i_1$.
Then we find some $\delta>0$ such that $\hat{x}_{i_1}- \delta>0$ and $\hat{x}_{\sigma(1)} +\delta < \min\{u_{\sigma(1)},m-M_1\}$.
Let $x \in \RR^n$ be defined as in Step 3. 
Then for all $i<\sigma(1)$ and all $i \ge i_1$ we have $\sum_{j =1}^i x_j = \sum_{j =1}^i \hat{x}_j$.
And for all $\sigma(1) \le i < i_1 $ we have 
\[
\sum_{j =1}^i x_j = \sum_{j =1}^i \hat{x}_j + \delta \quad \begin{cases}\quad \ge a_i \\ \quad = M_1 + \hat{x}_{\sigma(1)} + \delta < m \le b_i \end{cases} \,.
\]
Hence $x$ is feasible for~\eqref{eq:P}.
This leads to the same contradiction as at the end of Step 3.
\item In case $\hat{x}_{\sigma(1)} > M_2-M_1$ we show that $\hat{x}_{\sigma(1)} = 0$, from which the assertion follows.
Assume to the contrary that $\hat{x}_{\sigma(1)} > 0$.
Then we find some $\delta>0$ such that $\hat{x}_{\sigma(1)}- \delta>0$ and $\hat{x}_{\sigma(1)}- \delta > M_2-M_1$.
Now we define $x \in \RR^n$ by
\[
x_i:=\begin{cases} \hat{x}_{\sigma(1)}-\delta &, i=\sigma(1) \\ \hat{x}_i &, \mbox{otherwise} \end{cases} \,.
\]
Then for all $i<\sigma(1)$ we have $\sum_{j =1}^i x_j = \sum_{j =1}^i \hat{x}_j$.
And for all $i \ge \sigma(1)$ we have 
\[
\sum_{j =1}^i x_j = \sum_{j =1}^i \hat{x}_j - \delta \quad \begin{cases}\quad \le b_i \\ \quad \ge M_1 + \hat{x}_{\sigma(1)}-\delta > M_2 \ge a_i \end{cases} \,.
\]
Hence $x$ is feasible for~\eqref{eq:P} with
\[
\sum_{j \in I} c_j \cdot x_j = P_{min} - \delta \cdot c_{\sigma(1)} < P_{min} \,,
\]
which is a contradiction to the fact that $\hat{x}$ is a solution of~\eqref{eq:P}.
\end{enumerate}
\noindent \textbf{Step 5}: In case $c_{\sigma(1)}<0$ we prove that $\hat{x}_{\sigma(1)}=\min\{u_{\sigma(k)}, m-M_1\}$.
Assume to the contrary that $\hat{x}_{\sigma(1)} < \min\{u_{\sigma(1)},m-M_1\}$.
Then it follows from Step~3 that $\sum_{j=1}^{\sigma(1)-1} \hat{x}_j=M_1$.
If there is some index $i_1 > \sigma(1)$ such that $\hat{x}_{i_1}>0$, then we can argue as in the second case of Step 4 to arrive at a contradiction.
Hence we have $\hat{x}_i=0$ for all $i > \sigma(1)$.
Let $\delta>0$ be such that $\hat{x}_{\sigma(1)}+ \delta < \min\{u_{\sigma(1)},m-M_1\}$ and define $x \in \RR^n$ by
\[
x_i:=\begin{cases} \hat{x}_{\sigma(1)}+\delta &, i=\sigma(1) \\ \hat{x}_i &, \mbox{otherwise} \end{cases} \,.
\]
Then for all $i<\sigma(1)$ we have $\sum_{j =1}^i x_j = \sum_{j =1}^i \hat{x}_j$.
And for all $i \ge \sigma(1)$ we have 
\[
\sum_{j =1}^i x_j = \sum_{j =1}^i \hat{x}_j + \delta \quad \begin{cases}\quad \ge a_i \\ \quad = M_1 + \hat{x}_{\sigma(1)}+\delta < m \le b_i \end{cases} \,.
\]
Hence $x$ is feasible for~\eqref{eq:P} with
\[
\sum_{j \in I} c_j \cdot x_j = P_{min} + \delta \cdot c_{\sigma(1)} < P_{min} \,,
\]
which is a contradiction to the fact that $\hat{x}$ is a solution of~\eqref{eq:P}.

\noindent The above proof reveals that, when the $c_i$ are different and not equal to zero, the solution $\hat{x}$ of~\eqref{eq:P} is in fact unique and does not depend on the actual values $c_i$ but only on their order via $\sigma$ and whether they are positive or negative.
Note that Algorithm only uses this information, too.
It remains to show that Algorithm also computes a solution of~\eqref{eq:P} if some of the $c_i$ coincide or are equal to zero.
Fix a permutation $\sigma$ of the indices which corresponds to increasing $c_i$ values, i.e. $c_{\sigma(1)} \le \ldots \le c_{\sigma(n)}$, and define $c_{\sigma(j)}^\delta:=c_{\sigma(j)}+j \cdot \delta$ for $\delta>0$ and all $j \in I$.
Then we still have the same, but now strict, order $c_{\sigma(1)}^\delta < \ldots < c_{\sigma(n)}^\delta$, and for all $\delta>0$ small enough we also have $c_{\sigma(j)}^\delta \not=0$ for all $j \in I$ and positive (negative) values remain positive (negative).
Hence for all $\delta>0$ small enough all linear programs ($P^\delta$) with $c_i$ replaced by $c_i^\delta$ have the same unique solution $\hat{x}$, independent of $\delta$.
We show that $\hat{x}$ is also a solution of~\eqref{eq:P}.
Indeed, since~\eqref{eq:P} and (P$^\delta$) have the same constraints, for any solution $x$ of~\eqref{eq:P} we have
\[
P_{min} \le \sum_{j \in I} c_j \cdot \hat{x}_j = \sum_{j \in I} c_j^\delta \cdot \hat{x}_j - \delta \cdot \sum_{j \in I} j \cdot \hat{x}_j \le \sum_{j \in I} c_j^\delta \cdot x_j - \delta \cdot \sum_{j \in I} j \cdot \hat{x}_j \,.
\]
By letting $\delta \to 0$ we see that the right hand side converges to $\sum_{j \in I} c_j \cdot x_j=P_{min}$.
Hence we also have $\sum_{j \in I} c_j \cdot \hat{x}_j =P_{min}$, i.e. $\hat{x}$ is a solution of~\eqref{eq:P}.
Note that in this way a zero value $c_i=0$ is treated in the same way as a positive value, which is also reflected by the choice ``$\ge$'' in Line 5 of Algorithm.
Based on what we have said above there is no need to change the input values from $c_i$ to $c_i^\delta$, since the actual values do not matter.
Furthermore, we could as well use ``$>$'' in Line 5 of Algorithm to treat zero values in the same way as negative values, and a similar proof as above with $c_{\sigma(j)}^\delta:=c_{\sigma(j)}- \frac{\delta}{j}$ confirms that this modified Algorithm also computes a solution of~\eqref{eq:P}.

\subsection{Modifications} 
\label{subsec:appsubsec}

A slight modification of~\eqref{eq:P} with constraints like
\[a_i \le \sum_{j=1}^i e_i \cdot p_j \cdot x_j \le b_i, \ \mbox{$\forall \ i\in I$} \,,\]
and positive numbers \(e, p >0\) can be transformed in an equivalent linear program of the form (P), by setting
\[\tilde{x}_j:=p_j \cdot x_j\,, \quad \tilde{c}_j:=\frac{c_j}{p_j}\,, \quad \tilde{a_i}:=\frac{a_i}{e_i}\,, \quad \tilde{b}_i:=\frac{b_i}{e_i}\,, \quad \tilde{u}_j:=p_j \cdot u_j\]
Note that \(\tilde{c}^T \cdot \tilde{x}\) and \(c^T \cdot x\) are equal, but the order of the transformed values \(\tilde{c}_j\) may be different from the original order of the values \(c_j\). A special case are the constraints
\[a_i \le \sum_{j=1}^{i}{q^{i-j}\cdot x_j \le b_i}\,,\]
with \(q>0\), by setting \(e_i:=q^i\) and \(p_j:=q^{-j}\). These constraints can be used to model storage units including losses.

\section{Regression analyses}\label{app:regression}

In this section, we illustrate the regression analyses applied to interpolate the efficiencies in section~\ref{sec:sec2} as well as the cost functions in section~\ref{sec:sec3}.

To calculate hourly energy loss (or corresponding energy reduction) factors for different storage capacities, the energy content of the storage unit is assumed to decrease exponentially. For simplicity, the thermal storage medium is assumed to be water with a specific heat capacity of 4200 J/(kg$\cdot$K) and the maximum possible temperature difference of the medium is set to 60~Kelvin. This implies a specific energy capacity of 70~kWh/m$^3$. Additionally, the energy is assumed to be homogeneously distributed within the storage. As a data basis for efficiencies of different storage units, the capacities (in liters) and the corresponding daily losses (in kWh per day) were taken from specifications of state-of-the-art products of a manufacturer of heat storage units \cite[pp. 24-25]{vaillant2015}. For the regression, these capacities were transformed from liter into kWh using the specific energy capacity. The resulting data is listed in table \ref{table:regression}. As a consequence of the observed linear correlation between the logarithmic values of the capacities and the corresponding daily energy losses, an exponential regression analysis based on the least squares method is used to model hourly energy reduction factors for different capacities. The result of the regression yields to the approximated hourly energy reduction factors $q^{\text{hour}}_{\text{approx}}$, which are given by
\begin{equation}\label{eq:approx_losses}
q^{\text{hour}}_{\text{approx}} = \left( \frac{c - 0.2431954 \cdot c^{\,0.61876}}{c} \right) ^{\frac{1}{24}}\,,
\end{equation}
where $c$ is the capacity of the storage in kWh. This function is used for a capacity-specific consideration of losses in section \ref{sec:sec3}. The data of the supporting points from the product specifications and the resulting fitted hourly energy reduction factors are presented in table~\ref{table:regression}.

The procedure for the estimation of the annual costs for different charging and storage units used in sections \ref{subsec:subsec2} and \ref{subsec:subsec3} was as follows:~The costs for the charge units stems from a linear regression analysis applied on the net values of the data from~\cite[p. 26]{solarbayer2014}. The result is a linear function with a slope of about 7~[EUR/kW] and an offset of 130~EUR. For system optimization (i.e.\ minimization of \eqref{eq:nonlinear problem}) the constant offset of 130~EUR was not considered, since it does not have an impact on the optimal configuration, which is due to the linearity of the problem. Taking into account a typical lifetime of 15 years the annual costs for storage units, hence, were estimated at about 0.47~EUR/kW.

The cost-function for storage units is the result of the averaged net values for heat storage units, neglecting installation and maintenance costs, which leads to 0.8 -- 1.57~EUR/l \cite{konzepte2009, henning2013energiesystem}. For simplicity, the capacity costs are set to 1~EUR/l. In order to obtain annual costs per kWh, these costs were scaled with the storage unit specific energy capacity and an assumed lifetime of 15 years, resulting in costs of about 0.95~EUR/kWh.

Above considerations and data were used to define the specific cost function 
\eqref{eq:nonlinear problem example1} used in the examples from the more general function \eqref{eq:nonlinear problem}.

\begin{table}
\centering
\begin{tabular}{l|c|c|c|c|c|c|c}
\hline
product & 1 & 2 & 3 & 4 & 5 & 6 & 7 \\
\hline
\hline
capacity~[kWh] & 3.5 & 5.6 & 7.0 & 8.4 & 14.0 & 21.0 & 28.0 \\
daily losses~[kWh/24h] & 0.54 & 0.66 & 0.79 & 0.92 & 1.4 & 1.6 & 1.8 \\
hourly reduction factor &0.9930	&0.9948 & 0.9950 & 0.9952 & 0.9956 & 0.9967 & 0.9972\\
\hline
fitted factor $q^{\text{hour}}_{\text{approx}}$ & 0.9932 & 0.9944 & 0.9949 & 0.9952 & 0.9961 & 0.9967 & 0.9971 \\
\hline
\end{tabular}
\caption{Storage capacities and corresponding daily energy losses from a common manufacturer of heat storage units and the corresponding hourly energy reduction factors compared with the results from the exponential regression analysis, equation~\eqref{eq:approx_losses}.}
\label{table:regression}
\end{table}

\clearpage
\section{Python code}\label{app:python}
\begin{lstlisting}
#!/usr/bin/python
# -*- coding: utf-8 -*-

# (C) 2015 by Lars Siemer, NEXT ENERGY, EWE Research Center for 
# Energy Technology, 26129 Oldenburg, Germany. 
# Contact information: email: lars.siemer@next-energy.de

# Optimization algorithm to solve the following constrained
# linear program:
#                   min c*x, 
#   subject to:     0 <= x_i <= u_i 
#   and             a_i <= x_1+...+x_i <= b_i, for all i=1,...,n.
# This code evolved in connection with the article
#   "Cost-Optimal Operation of Energy Storage Units: 
#        Benefits of a Problem-Specific Approach"
#==================================================================

import numpy

#==================================================================
# definition of the algorithm
#==================================================================

def Algorithm(a,b,c,u):
    sigma=numpy.argsort(c)
    dim=numpy.size(c)
    x=numpy.zeros(dim)
    for k in range(0,dim):
        i=sigma[k]
        M_1=numpy.max(numpy.insert(a[:i],0,0))
        M_2=numpy.max(numpy.insert(a[i:],0,0))
        m=numpy.min(b[i:])
        if c[i]<0:
            x[i]=numpy.min([u[i],m-M_1])
        else:
            x[i]=numpy.min([numpy.max([0,M_2-M_1]),\
                 numpy.min([u[i],m-M_1])])
        for l in range(i,dim):
            a[l]=a[l]-x[i]
            b[l]=b[l]-x[i]
    return x
\end{lstlisting}

\begin{thebibliography}{10}

\bibitem{georgilakis2008technical}
P.~S. Georgilakis.
\newblock Technical challenges associated with the integration of wind power
  into power systems.
\newblock {\em Renewable and Sustainable Energy Reviews}, 12(3):852--863, 2008.

\bibitem{heide2010seasonal}
D.~Heide, L.~Von~Bremen, M.~Greiner, C.~Hoffmann, M.~Speckmann, and
  S.~Bofinger.
\newblock Seasonal optimal mix of wind and solar power in a future, highly
  renewable {Europe}.
\newblock {\em Renewable Energy}, 35(11):2483--2489, 2010.

\bibitem{rasmussen2012storage}
M.~G. Rasmussen, G.~B. Andresen, and M.~Greiner.
\newblock Storage and balancing synergies in a fully or highly renewable
  pan-{European} power system.
\newblock {\em Energy Policy}, 51:642--651, 2012.

\bibitem{weitemeyer2015integration}
S.~Weitemeyer, D.~Kleinhans, T.~Vogt, and C.~Agert.
\newblock Integration of renewable energy sources in future power systems: The
  role of storage.
\newblock {\em Renewable Energy}, 75:14--20, 2015.

\bibitem{potenzial2014}
D.~B{\"o}ttger, M.~G{\"o}tz, N.~Lehr, H.~Kondziella, and T.~Bruckner.
\newblock Potential of the power-to-heat technology in district heating grids
  in {Germany}.
\newblock {\em Energy Procedia}, 46:246--253, 2014.

\bibitem{sterner2014energiespeicher}
Michael Sterner and Ingo Stadler.
\newblock {\em Energiespeicher-Bedarf, Technologien, Integration}.
\newblock Springer-Verlag, 2014.

\bibitem{denholm2010role}
P.~Denholm, E.~Ela, B.~Kirby, and M.~Milligan.
\newblock The role of energy storage with renewable electricity generation.
\newblock {\em National Renewable Energy Laboratory (NREL)}, 2010.

\bibitem{pedersen2011}
T.~S. Pedersen, P.~Andersen, H.~L. Nielsen, K. M.and~Starmose, and P.~D.
  Pedersen.
\newblock Using heat pump energy storages in the power grid.
\newblock In {\em Control Applications (CCA), 2011 IEEE International
  Conference on}, pages 1106--1111. IEEE, 2011.

\bibitem{eiselt2012dezentrale}
J.~Eiselt.
\newblock {\em {Dezentrale Energiewende: Chancen und Herausforderungen}}.
\newblock Springer-Verlag, 2012.

\bibitem{stat2014energie}
Statistisches Bundesamt.
\newblock
  {Umweltnutzung
  und Wirtschaft - Tabellen zu den Umwelt{\"o}konomischen Gesamtrechnungen Teil
  2: Energie}.
\newblock 2014.

\bibitem{nocedal2006numerical}
J.~Nocedal and S.~J. Wright.
\newblock {\em Numerical Optimization}.
\newblock Springer, 2 edition, 2006.

\bibitem{von1981heat}
H.~L. Von~Cube and F.~Steimle.
\newblock {\em Heat pump technology}.
\newblock Butterworth-Heinemann Ltd, 1981.

\bibitem{perko2011calculation}
J.~Perko, V.~Dugec, D.~Topic, D.~Sljivac, and Z.~Kovac.
\newblock Calculation and design of the heat pumps.
\newblock In {\em Energetics (IYCE), Proceedings of the 2011 3rd International
  Youth Conference on}, pages 1--7. IEEE, 2011.

\bibitem{dincer2002thermal}
I.~Dincer and M.~Rosen.
\newblock {\em Thermal energy storage: systems and applications}.
\newblock John Wiley \& Sons, 2002.

\bibitem{schroder1981latent}
J.~Schr{\"o}der and K.~Gawron.
\newblock Latent heat storage.
\newblock {\em International Journal of Energy Research}, 5(2):103--109, 1981.

\bibitem{zensus2011b}
Statistische~{\"A}mter des Bundes und~der L{\"a}nder.
\newblock
  {Zensus
  2011 - Geb{\"a}ude- und Wohnungsbestand in Deutschland -- Erste Ergebnisse
  der Geb{\"a}ude- und Wohnungsz{\"a}hlung 2011}.
\newblock 2014.

\bibitem{zensus2011a}
Statistische~{\"A}mter des Bundes und~der L{\"a}nder.
\newblock {Zensus 2011 - Haushalte und Familien -- Ergebnisse des Zensus am 9. Mai 2011}.
\newblock 2014.

\bibitem{entwicklung2010}
Statistische Bundesamt.
\newblock
  {Bev{\"o}lkerung
  und Erwerbst{\"a}tigkeit -- Entwicklung der Privathaushalte bis 2030~--
  Ergebnisse der Haushaltsvorausberechnung}.
\newblock 2011.

\bibitem{bigalke2012dena}
U.~Bigalke, H.~Discher, H.~Lukas, Y.~Zeng, K.~Bensmann, and C.~Stolte.
\newblock {Der dena Geb{\"a}udereport 2012}.

\bibitem{EnEV2009}
Bundesministerium der Justiz.
\newblock
  {Bundesgesetzblatt
  -- Verordnung zur {\"A}nderung der Energieeinsparverordnung}.
\newblock {\em Bundesgesetzblatt Teil 1}, (23):954--992, 2009.

\bibitem{directive2010}
EU~Directive.
\newblock
  {Directive
  2010/30/EU of the European Parliament and of the council of 19 May 2010 on
  the indication by labelling and standard product information of the
  consumption of energy and other resources by energy-related products}.
\newblock {\em European Parliament and the Council of the European Union},
  2010.

\bibitem{vaillant2015}
Vaillant.
\newblock {\em
  {Elektro-Warmwasser:
  VED, VEH, VEN, VEK}}, 2015.

\bibitem{solarbayer2014}
Solarbayer GmbH.
\newblock {\em
  {Preisliste
  - Bruttopreisliste 2014/2015}}, 2014.

\bibitem{konzepte2009}
C.~Wilhelms, K.~Za\ss, K.~Vajen, and U.~Jordan.
\newblock {\em Neue Konzepte f{\"u}r Warmwasserspeicher bis
  50m{$^3$}~--~Markt{\"u}bersicht, Kosten, Anwendungsgebiete}, 2009.
\newblock VDI--Bericht Nr. 2074.

\bibitem{henning2013energiesystem}
H.~M. Henning and A.~Palzer.
\newblock
  {Energiesystem Deutschland 2050. Sektor-und Energietr{\"a}ger{\"u}bergreifende, modellbasierte, ganzheitliche Untersuchung zur langfristigen Reduktion energiebedingter CO{$_2$}-Emissionen durch Energieeffizienz und den Einsatz Erneuerbarer Energien}.
\newblock {\em FhG-ISE, Studie im Rahmen eines BMWi-Projektes, Freiburg}, 2013.

\end{thebibliography}
\end{document}